\documentclass{ifacconf}%
\usepackage{color}
\usepackage{natbib}
\usepackage{epsfig}
\usepackage{amsmath}
\usepackage{amsfonts}
\usepackage{subfigure}
\usepackage{graphicx}
\usepackage{amssymb}%
\newcommand{\bx}{\mathrm{\mathbf{x}}}

\newcommand{\bu}{\mathrm{\mathbf{u}}}

\newcommand{\btheta}{\boldsymbol\theta}

\newcounter{asscount}
\newenvironment{assumption}[1][]{\refstepcounter{asscount}\par\noindent   \textbf{Assumption~\theasscount. #1 } \rmfamily  }{\par}
\newcounter{remcount}

\newcounter{lemcount}

\newcounter{themcount}
\newenvironment{theorem}[1][]{\refstepcounter{themcount}\par   \textbf{Theorem~\thethemcount. #1 } \rmfamily }{\par}
\begin{document}
\begin{frontmatter}
\title{A General Framework for Modeling and Online Optimization of Stochastic Hybrid Systems\thanksref{footnoteinfo}}
\thanks[footnoteinfo]{The authors' work is supported in part by NSF under
Grant EFRI-0735974, by AFOSR under grant FA9550-09-1-0095, by DOE
under grant DE-FG52-06NA27490, by ONR under grant
N00014-09-1-1051, and by ARO under grant W911NF-11-1-0227.}
\author[First]{A. Kebarighotbi}
\author[First]{C.~G. Cassandras}
\address[First]{Division of Systems Engineering, Center for Information and Systems Engineering,  Boston University, Brookline, MA 02446 \texttt{alik@bu.edu,cgc@bu.edu}}
\begin{abstract}
We extend the definition of a Stochastic Hybrid Automaton (SHA) to
overcome limitations that make it difficult to use for on-line
control. Since guard sets do not specify the exact event causing a
transition, we introduce a clock structure (borrowed from timed
automata), timer states, and guard functions that disambiguate how
transitions occur. In the modified SHA, we formally show that
every transition is associated with an explicit element of an
underlying event set. This also makes it possible to uniformly
treat all events observed on a sample path of a stochastic hybrid
system and generalize the performance sensitivity estimators
derived through Infinitesimal Perturbation Analysis (IPA). We
eliminate the need for a  case-by-case treatment of different
event types and provide a unified set of matrix IPA equations. We
illustrate our approach by revisiting an optimization problem for
single node finite-capacity stochastic flow systems to obtain
performance sensitivity estimates in this new setting.
\end{abstract}
\end{frontmatter}

\section{Introduction}

A Stochastic Hybrid System (SHS) consists of both time-driven and
event-driven components. Its stochastic features may include
random transition times and external stochastic inputs or
disturbances. The modeling and optimization of these systems is
quite challenging and many models have been proposed, some
capturing randomness through probabilistic resets when reset
functions are distributions, through spontaneous transitions at
random times \citep{BujLyg06},\citep{Hespanha04}, Stochastic
Differential Equations (SDE)
\citep{GhoshArapMarc93},\citep{GhoshArapMarc97}, or using
Stochastic Flow Models (SFM) \citep{CassWarMelSunPan02} with the
aim of describing stochastic continuous dynamics.

Optimizing the performance of SHS poses additional challenges and
most approaches rely on approximations and/or using
computationally taxing methods. For example, \citep{BujLyg04} and
\citep{Koutsoukos05} resort to dynamic programming techniques. The
inherent computational complexity of these approaches makes them
unsuitable for on-line optimization. However, in the case of
parametric optimization, application of \emph{Infinitesimal
Perturbation Analysis} (IPA) \citep{CassWarMelSunPan02} to SHS has
been very successful in on-line applications. Using IPA, one
generates sensitivity estimates of a performance objective with
respect to a control parameter vector based on readily available
data observed from a single sample path of the SHS. Along this
line, SFMs provide the most common framework for applying IPA in
the SHS setting and have their root in making abstractions of
complex \emph{Discrete Event Systems} (DES), where the event rates
are treated as stochastic processes of arbitrary generality except
for mild technical assumptions. A fundamental property of IPA is
that the derivative estimates obtained are independent of the
probability laws of the stochastic processes involved. Thus, they
can be easily obtained and, unlike most other techniques, they can
be implemented in on-line algorithms.

In this paper, we aim at extending \emph{Stochastic Hybrid
Automata} (SHA) \citep{SHSbook} and create a framework within
which IPA becomes straightforward and applicable to arbitrary SHS.
A SHA, specifies discrete states (or \emph{modes}) where the state
$\mathrm{\mathbf{x}}$ evolves according to a continuous vector
field until an event triggers a mode transition. The transitions
are described by \emph{guards} and \emph{invariants} as well as
clock structures \citep{BujLyg06} borrowed from \emph{Stochastic
Timed Automata} (STA) \citep{glynn89},\citep{cassbook}. When
$\mathrm{\mathbf{x}}$ reaches a guard set, a transition becomes
\emph{enabled} but not triggered. On the other hand, if
$\mathrm{\mathbf{x}}$ exits the boundary defined by the invariant
set in a mode or if a spontaneous event occurs, the transition
must trigger. This setting has the following drawbacks: $(i)$ The
clock structure (normally part of the system input) is not
incorporated in the definition of guards and invariant conditions.
As a result, spontaneous transitions have to be treated
differently. $(ii)$ The guard set does not specify the exact event
causing a transition and we cannot, therefore, differentiate
between an event whose occurrence time may depend on some control
parameter $\boldsymbol{\theta}$ and another that does not. This is
a crucial point in IPA, as it directly affects how performance
derivative estimates evolve in time. As described in
\citep{CassWarPanYao09}, when applied to SFMs, IPA uses a
classification of events into different types (exogenous,
endogenous, and induced) to extract this information.

Here, we seek an enriched model which explicitly specifies the
event triggering a transition and, at the same time, creates a
unified treatment for all events, i.e., all IPA equations are
common regardless of event type. We achieve this by introducing
state variables representing \emph{timers} and treating a clock
structure as an input to the system with the mode invariants
generally dependent on both. We formalize the definition of an
event by associating it to a \emph{guard function} replacing the
notion of the guard set. This removes the ambiguity caused by an
enabled, but not triggered, event, as well as the need for
treating spontaneous transitions differently. A byproduct of this
unified treatment is the development of a matrix notation for the
IPA equations, which makes the treatment of complex systems with
multiple states and events a straightforward application of these
equations. We verify this process by applying it to a single node
queueing system previously solved using the SFM framework
\citep{CassWarMelSunPan02}.

The paper is organized as follows: In Section 2, we present the
general SHA which includes all the features previously handled by
SFMs. Utilizing the resulting model, in Section 3 we develop a
matrix notation for IPA which simplifies the derivation of
sensitivity estimators in Section 4. We verify our results in
Section 5 by applying the proposed technique to a single node
finite capacity buffer system. We conclude with Section 6.

\section{General Optimization Model}

Let us consider a \emph{Stochastic Hybrid Automaton} (SHA) as
defined in \citep{cassbook} with only slight modifications and
parameterized by the
vector $\boldsymbol{\theta}=(\theta_{1},\theta_{2},\ldots,\theta_{N_{\theta}%
})\in\Theta$ as
\[
G=(\mathcal{Q},\mathcal{X},\mathcal{U},\Theta,\mathcal{E}%
,\mathrm{\mathbf{f}},\phi,Inv,guard,\mathrm{\mathbf{r}},(q_{0}%
,\mathrm{\mathbf{x}}_{0}))
\]
\noindent where\newline\noindent$\bullet\ $ $\mathcal{Q}\subset\mathbb{Z}^{+}$
is the countable set of discrete states or modes;\newline\noindent$\bullet\ $
$\mathcal{X}(\boldsymbol{\theta})\subset\mathbb{R}^{N_{x}}$ is the admissible
continuous state space for any $\boldsymbol{\theta}\in\Theta$;\newline%
\noindent$\bullet\ $ $\mathcal{U}(\boldsymbol{\theta})\subset\mathbb{R}%
^{N_{u}}$ is the set of inputs (possibly disturbances or clock variables) for
any $\boldsymbol{\theta}\in\Theta$;\newline\noindent$\bullet\ $ $\Theta
\subset\mathbb{R}^{N_{\theta}}$ is the set of admissible control
parameters;\newline\noindent$\bullet\ $ $\mathcal{E}$ is a countable event set
$\mathcal{E}=\{E_{i},i=1,2,\ldots,N_{e}\}$;\newline\noindent$\bullet\ $
$\mathrm{\mathbf{f}}$ is a continuous vector field, $\mathrm{\mathbf{f}%
}:\mathcal{Q}\times\mathcal{X}({\boldsymbol{\theta}})\times\mathcal{U}%
({\boldsymbol{\theta}})\times\Theta\mapsto\mathcal{X}({\boldsymbol{\theta}}%
)$;\newline\noindent$\bullet\ $ $\phi$ is a discrete transition function
$\phi:\mathcal{Q}\times\mathcal{X}({\boldsymbol{\theta}})\times\mathcal{U}%
(\boldsymbol{\theta})\times{\mathcal{E}}\mapsto\mathcal{Q}$;\newline%
\noindent$\bullet\ $ $Inv$ is a set defining an invariant condition such that
$Inv\subseteq\mathcal{Q}\times\mathcal{X}({\boldsymbol{\theta}})\times
\mathcal{U}(\boldsymbol{\theta})\times\Theta$;\newline\noindent$\bullet$
$guard$ is a set defining a guard condition, $guard\subseteq\mathcal{Q}%
\times\mathcal{Q}\times\mathcal{X}(\boldsymbol{\theta})\times\mathcal{U}%
(\boldsymbol{\theta})\times\Theta$;\newline\noindent$\bullet\ $
$\mathrm{\mathbf{r}}$ is a reset function, $\mathrm{\mathbf{r}}:\mathcal{Q}%
\times\mathcal{Q}\times\mathcal{E}\times\mathcal{X}({\boldsymbol{\theta}%
})\times\mathcal{U}(\boldsymbol{\theta})\times\Theta\mapsto\mathcal{X}%
({\boldsymbol{\theta}})$;\newline\noindent$\bullet\ $ $(q_{0}%
,\mathrm{\mathbf{x}}_{0})$ is the initial state.

Note that the input $u\in\mathcal{U}(\btheta)$ can be a vector of
random processes which are all defined on a common probability
space $(\Omega,\mathcal{F},P)$. Also, observe that invariants and
guards are sets which do not specify the events (hence, precise
times) causing violation or adherence to their set conditions.
Thus, we cannot differentiate between a transition that depends on
$\boldsymbol{\theta}$ and one that does not. This prevents us from
properly estimating the effect of $\boldsymbol{\theta}$ on the
system behavior. In particular, if
$(\mathrm{\mathbf{x}},\mathrm{\mathbf{u}})\in guard(q,q^{\prime
})$ for some $q,q^{\prime}\in\mathcal{Q}$, a transition to $q^{\prime}%
\in\mathcal{Q}$ can occur either through a policy that uniquely specifies some
$(\mathrm{\mathbf{x}},\mathrm{\mathbf{u}},\boldsymbol{\theta})$ causing the
transition or at some random time while $(\mathrm{\mathbf{x}}%
,\mathrm{\mathbf{u}},\boldsymbol{\theta})$ remains in the guard
set. This is one of the issues we focus on in what follows.

We allow the parameter vector $\boldsymbol{\theta}$ to affect the system not
only through the vector field, reset conditions, guards, and invariants, but
also its structure through $\mathcal{X}({\boldsymbol{\theta}})$ and
$\mathcal{U}({\boldsymbol{\theta}})$, e.g., the parameters can appear in the
state and input constraints. We remove $\boldsymbol{\theta}$ from the
arguments whenever it does not cause any confusion and simplifies notation.
Defining $\mathrm{\mathbf{x}}(t,\boldsymbol{\theta})$ and $\mathrm{\mathbf{u}%
}(t,\boldsymbol{\theta})$ as the state and input vectors, we introduce the
following assumptions:

\begin{assumption}
With probability 1, for any $\mathrm{\mathbf{x}}\in\mathcal{X}%
,q\in\mathcal{Q}, \mathrm{\mathbf{u}}\in\mathcal{U}$, and
$\boldsymbol{\theta}\in\Theta$,
$\|\mathrm{\mathbf{f}}(q,\bx,\bu,\btheta)\|_{\infty }<\infty$
where $\|\cdot\|_{\infty}$ is the $L_{\infty}$
norm.\label{assump:finitef}
\end{assumption}

\begin{assumption}
With probability 1, no two events can occur at the same time unless one causes
the occurrence of the other.\label{assump:nosimul}
\end{assumption}

Assumption \ref{assump:finitef} ensures that $\mathrm{\mathbf{x}%
}(t,\boldsymbol{\theta})$ remains smooth inside a mode as is embedded in the
definition of the SHA above. Assumption \ref{assump:nosimul} rules out the
pathological case of having two independent events happening at the same time.

\begin{figure}[th]
\centering
\includegraphics[scale =.3]{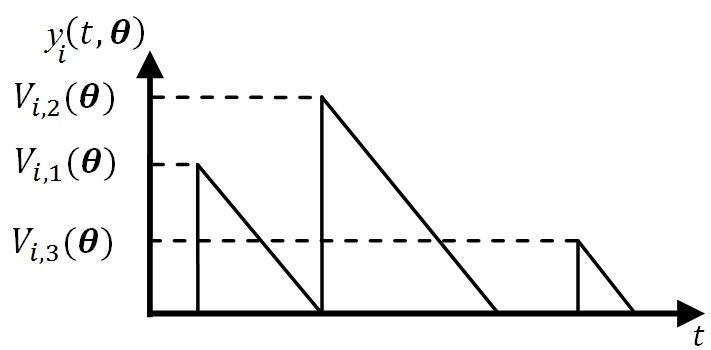}\caption{{\small A
timer realization based on a given clock sequence }{\footnotesize {$\{V_{i,n}%
\}$}}}%
\label{fig:timerinput}%
\end{figure}

Borrowing the concept of clock structure from \emph{Stochastic Timed Automata
}as defined in \citep{cassbook}, we associate an event $E_{i}\in\mathcal{E}$
with a a sequence $\{V_{i,1}(\boldsymbol{\theta}),V_{i,2}(\boldsymbol{\theta
}),\ldots\}$ where $V_{i,n}(\boldsymbol{\theta})$ is the $n$th (generally
random) lifetime of $E_{i}$, i.e., if this event becomes feasible for the
$n$th time at some time $t$, then its next occurrence is at $t+V_{i,n}%
(\boldsymbol{\theta})$. Obviously, not all events in $\mathcal{E}$
are defined to occur in this fashion, but if they are, we define a
\emph{timer} as a state
variable, say $y_{i}$, so that it is initialized to $y_{i}(t)=V_{i,n}%
(\boldsymbol{\theta})$ if $E_{i}$ becomes feasible for the $n$th time at time
$t$. Subsequently, the timer dynamics are given by $\dot{y}_{i}=-1$ until the
timer runs off, i.e., $y_{i}(t+V_{i,n}(\boldsymbol{\theta}))=0$. Figure
\ref{fig:timerinput} shows an example of a timer state as it evolves according
to the supplied event lifetimes. We assume that $V_{i,n}(\boldsymbol{\theta})$ is differentiable
with respect to $\boldsymbol{\theta}$ for all $n$.

\textbf{The concept of event. }In the SHA as defined above, $\mathcal{E}$ is
simply a set of labels. We provide more structure to an event by assigning to
each $E_{i}\in\mathcal{E}$, a \emph{guard function} $g_{i}:\mathbb{R}%
^{+}\times\mathcal{X}\times\mathcal{U}\times\Theta\mapsto\mathbb{R}$
which is not null (i.e. $g_{i}\not \equiv 0$) and is assumed to be
differentiable almost everywhere on its domain. We are interested
in a sample path of a SHA $G$ on some interval $[0,T]$ where we
let $\tau_{k}(\boldsymbol{\theta})$ be the time when the $k$th
transition fires and set $0=\tau_{0}\leq\tau
_{1}(\boldsymbol{\theta})\leq\ldots\leq\tau_{K}(\boldsymbol{\theta})\leq
\tau_{K+1}=T$. We then define an event $E_{i}$ as occurring at
time $\tau _{k}(\boldsymbol{\theta})$ if
\[
\tau_{k}=\inf\{t\geq\tau_{k-1}:g_{i}(t,\mathrm{\mathbf{x}},\mathrm{\mathbf{u}%
},\boldsymbol{\theta})=0\}
\]
that is, the event satisfies the condition $g_{i}(t,\mathrm{\mathbf{x}%
},\mathrm{\mathbf{u}},\boldsymbol{\theta})=0$ which was not being
satisfied over
$(\tau_{k-1}(\boldsymbol{\theta}),\tau_{k}(\boldsymbol{\theta}))$.
The following theorem shows that using guard functions, we can
associate every transition with an event occurring at the
transition time.
\begin{theorem}
For every STA $G$ with $\mathcal{E}$, $\mathcal{X}$ and
$\mathcal{Q}$, there exists another STA $\tilde{G}$ with event set
$\tilde{\mathcal{E}}$, continuous state space
$\tilde{\mathcal{X}}$ and discrete state space
$\tilde{\mathcal{Q}}$ such that $\tilde{\mathcal{Q}}=\mathcal{Q}$
and every transition $(q,q^{\prime})$ in $G$ is associated with an
event $e\in\tilde{\mathcal{E}}$ at the transition
time.\label{thm:existence}
\end{theorem}

\textit{Proof:} A transition $(q,q^{\prime})$ is dictated by the
transition function
$\phi(q,\mathrm{\mathbf{x}},\mathrm{\mathbf{u}},e)$ such that
$\phi(q,\mathrm{\mathbf{x}},\mathrm{\mathbf{u}},e)=q^{\prime}$ for
some
$\mathrm{\mathbf{x}}\in\mathcal{X}(\boldsymbol{\theta}),$ $\mathrm{\mathbf{u}%
}\in\mathcal{U}(\boldsymbol{\theta}),$ $e\in\mathcal{E}$. If
$\phi(q,\bx,\bu,E_i)=\phi (q,E_{i})=q^{\prime}$ for some
$E_{i}\in\mathcal{E}$, i.e., the transition
$(q,q^{\prime})$ is independent of $\mathrm{\mathbf{x}}\in\mathcal{X}%
(\boldsymbol{\theta}),$
$\mathrm{\mathbf{u}}\in\mathcal{U}(\boldsymbol{\theta })$, the
proof is complete. In this case, we can always augment $\bx \in
\mathcal{X}$ to $\tilde{\bx}=(\bx,y_i)\in \tilde{\mathcal{X}}$
where $y_i$ is a timer state variable  capturing lifetimes of
event $E_i\in\mathcal{E}$ and
associate $E_i$ with  guard function $g_i(t,\tilde\bx,\bu,\btheta)=y_i$. If $\phi(q,\bx,\bu,e)=\phi(q,\mathrm{\mathbf{x}},\mathrm{\mathbf{u}%
})=q^{\prime}$, i.e., the transition $(q,q^{\prime})$ depends on
$(\mathrm{\mathbf{x}},\mathrm{\mathbf{u}},\boldsymbol{\theta})$,
then it is
either a result of violating $Inv(q)$ or it occurs while $(\mathrm{\mathbf{x}%
},\mathrm{\mathbf{u}},\boldsymbol{\theta})\in
guard(q,q^{\prime})$. In the
former case, we can define some $E_{i}$ such that $g_{i}(t,\mathrm{\mathbf{x}%
},\mathrm{\mathbf{u}},\boldsymbol{\theta})=0$ is the condition
that determines the occurrence time of $E_{i}$. This is because
$Inv(q)$ can be violated in two ways: $(a)$ directly, due to an
occurrence of $E_{i}$ meaning $(\bx,\bu)$ is on the boundary of
$Inv(q)$ at the transition time; $(b)$ indirectly, due to a reset
condition which is the result of a previous transition
$\phi(q^{\prime },\mathrm{\mathbf{x}},\mathrm{\mathbf{u}},e)=q$,
where it is possible that $q^{\prime}=q$ (a self-loop transition).
That is, the reset condition is such
that $\mathrm{\mathbf{r}}(q^{\prime},q,e,\mathrm{\mathbf{x}}%
,\mathrm{\mathbf{u}},\boldsymbol{\theta})\notin Inv(q)$. In case $(a)$, $g_{i}(\cdot)$ is such that $g_i(t,\mathrm{\mathbf{x}%
},\mathrm{\mathbf{u}},\boldsymbol{\theta})=0$ is part of the
boundary of $Inv(q)$ including $(\bx,\bu)$. In case $(b)$, the
transition can only occur as $(i)$ a result of some
$e\in\mathcal{E}$ (completing the proof); $(ii)$ the violation of
$Inv(q^{\prime})$; or $(iii)$ while
$(\mathrm{\mathbf{x}},\mathrm{\mathbf{u}})\in
guard(q^{\prime},q)$. Since we have already considered the first
two cases, we only need to check case $(iii)$, including
$q^{\prime}=q$. This case can occur either
through $(A)$ a policy equivalent to a condition $g(t,\mathrm{\mathbf{x}%
},\mathrm{\mathbf{u}},\boldsymbol{\theta})=0$; or $(B)$ after some
random
time. In the former case, we can define some $E_{i}$ such that $g_{i}%
(t,\mathrm{\mathbf{x}},\mathrm{\mathbf{u}},\boldsymbol{\theta}%
)=g(t,\mathrm{\mathbf{x}},\mathrm{\mathbf{u}},\boldsymbol{\theta})=0$
and include $E_{i}$ $\in\tilde{\mathcal{E}}$. In case $(B)$, let
$\tau=\inf\{t\geq
\tau_{k-1}:(\mathrm{\mathbf{x}}(t),\mathrm{\mathbf{u}}(t))\in
guard(q^{\prime },q)\}$ be the time that $(\bx,\bu)$ enters
$guard(q^{\prime},q)$. We can associate a self-loop transition
$(q^\prime,q^\prime)$ at $\tau$ which is caused by some event
$E_{i}$ with guard
function $g_{i}(\cdot)$ such that $g_{i}(\tau,\mathrm{\mathbf{x}%
},\mathrm{\mathbf{u}},\boldsymbol{\theta})=0$. Note that
$(\bx,\bu)$ satisfying this condition forms part of the boundary
of $guard(q,q^\prime)$. We define a reset condition for this
transition such that a timer with state $y$ gets initialized with
a random value $V_i(\btheta)$ such that
$\tau+V_{i}(\boldsymbol{\theta})$ is the time of transition
$(q^\prime,q)$. We can then include $y$ in $\tilde{\bx}$ and
define the event at $\tau+V_i(\theta)$ as
$E_{j}\in\tilde{\mathcal{E}}$ and associate with it a guard
function $g_{j}(t,\tilde\bx,\bu,\boldsymbol{\theta})=y_{j}$. Since
such events $E_{i},E_{j}$ can always be defined, the proof is
complete.\hfill$\rule{2mm}{2mm}$

In the above proof, we turn the reader's attention to how timers
and guard functions remove the need for guard sets. Also, in the
case of a chain of simultaneous transitions, they identify an
event whose guard function determines when the transitions occur.

\begin{figure}[th]
\centering
\includegraphics[scale =.30]{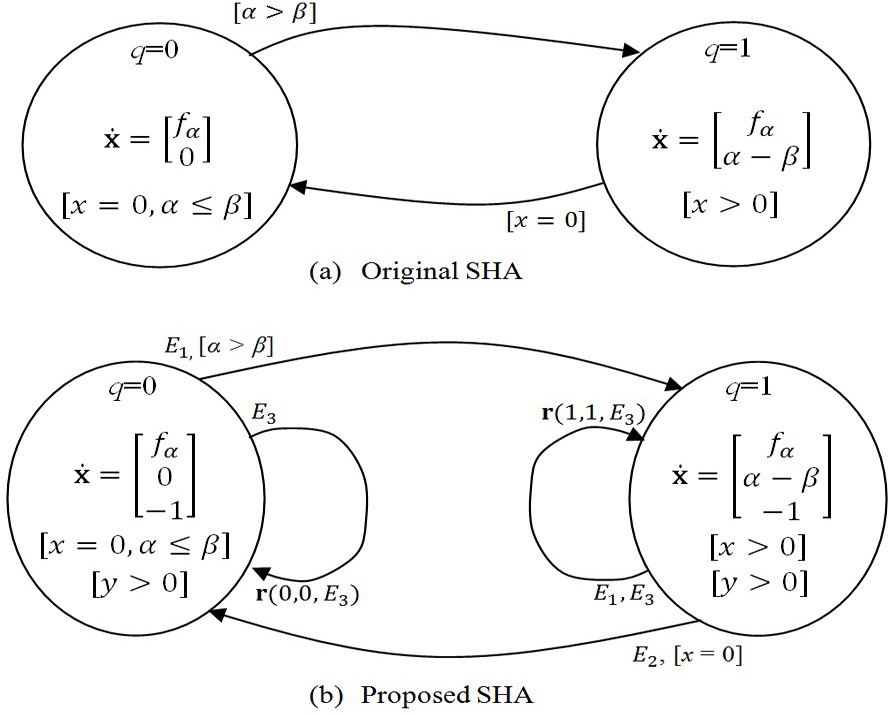}\caption{{\small A
simple fluid buffer SHA: Contrasting two approaches.}}%
\label{fig:twoSHAs}%
\end{figure}

\textbf{Example}. To illustrate this framework based on events
associated with guard functions, consider the example of a SHA as
shown in Fig. \ref{fig:twoSHAs}(a). This models a simple flow
system with a buffer whose content is $x(t,\theta)\geq0$. The
dynamics of the content are $\dot {x}(t,\theta)=0$ when $q=0$
(empty buffer) and $\dot{x}(t,\theta )=\alpha(t,\theta)-\beta$,
otherwise. Here, $\beta\geq0$ is a fixed outflow rate and
$\{\alpha(t,\theta)\}$ is a piecewise differentiable random
process whose behavior depends on $\theta$ via a continuous vector
field $f_{\alpha }(t,\theta)$. We allow for discontinuous jumps in
the value of $\alpha (t,\theta)$ at random points in time modeled
through events that occur at time instants
$V_{1},V_{1}+V_{2},\ldots$ using a timer state variable $y(t)$
re-initialized to $V_{n+1}$ after the $n$th time that $y(t)=0$.
The result of such an event at time $t$ is a new value
$\alpha(t^{+})=A_{n+1}$ where this jump process is independent of
$\theta$. In mode $q=0$, the invariant condition
$\alpha(t,\theta)\leq\beta$ is required to ensure that the buffer
remains empty.

The state of the new SHA, shown in Fig. \ref{fig:twoSHAs}(b), is
denoted by $\mathrm{\mathbf{x}}=(\alpha,x,y)$ where note that
$\dot{y}(t)=-1$. The event set is
$\mathcal{E}=\{E_{1},E_{2},E_{3}\}$ with
$g_{1}(t,\tilde\bx,\bu,\btheta)=\alpha(t,\theta)-\beta$,
$g_{2}(t,\tilde\bx,\bu,\btheta)=x(t,\theta)$, and
$g_{3}(t,\tilde\bx,\bu,\btheta)=y(t,\theta)$. In addition,
we define the reset condition $\mathrm{\mathbf{r}}(0,0,E_{3}%
)=\mathrm{\mathbf{r}}(1,1,,E_{3})=(A_{n+1},x,V_{n+1})$ whenever $E_{3}$ occurs
for the $n$th time, treating $\{A_{n}\}$ and $\{V_{n}\}$, $n=1,2,\ldots$, as
input processes.

When $q=0$, two events are possible: $(i)$ If $E_{1}$ occurs at $t=\tau_{k}$
we have $\alpha(\tau_{k},\theta)-\beta=0$ and a transition to $q=1$ occurs
since the condition $\alpha(t,\theta)-\beta<0$ must have held at $\tau_{k}%
^{-}$. $(ii)$ If $E_{3}$ occurs at $t=\tau_{k}$ we have $y(\tau_{k})=0$ and a
self-loop transition results. By the reset condition, $\alpha(\tau_{k}%
^{+})=A_{n+1}$ for some $A_{n+1}$, assuming this was the $n$th occurrence of
this event. If $A_{n+1}>\beta$, then immediately a transition to $q=1$ occurs.
Observe that even though the condition of this transition is $\alpha(\tau
_{k}^{+})>\beta$, the transition is still due to event $E_{3}.$

When $q=1$, all three events are possible, but $E_{1}$,$E_{3}$
cause self loops. If $E_{2}$ occurs at some $t=\tau_{k}$, then
$x(\tau_{k})=0$. On the other hand, suppose that a transition to
$q=1$ occurred because of $E_{1}$ at $t=\tau_{k}$, i.e.,
$\alpha(\tau_{k},\theta)-\beta=0$. It is possible that
$\alpha(t,\theta)-\beta=0$ for
$t\in\lbrack\tau_{k},\tau_{k}+\epsilon]$, $\epsilon>0$. In this
case, $\dot{x}=\alpha(t,\theta)-\beta=0$ and $x(\tau _{k}^{+})=0$.
This violates the invariant condition $[x>0]$ at $q=1$, causing an
immediate return to $q=0$. Similarly, if
$\alpha(\tau_{k},\theta)-\beta=0$ but
$\alpha(\tau_{k}^{+},\theta)-\beta<0$, the invariant condition at
$q=1$ is violated and there is an immediate return to $q=0$. All
these can be summarized in the transition functions below:
\begin{align*}
\phi(0,\mathrm{\mathbf{x}},\bu,e) &  =\left\{
\begin{array}
[c]{ll}%
1 & \mathrm{\ if\ }e=E_{1}\text{ or }\alpha(t,\theta)>\beta\\
0 & \mathrm{\ otherwise}%
\end{array}
\right.  \\
\phi(1,\mathrm{\mathbf{x}},\bu,e) &  =\left\{
\begin{array}
[c]{ll}%
0 & \mathrm{\ if\ }e=E_{2}\text{ or }x(t,\theta)=0\\
1 & \mathrm{\ otherwise}%
\end{array}
\right.
\end{align*}
where $\bu=\alpha$.

It is important to note that the conditions $[\alpha>\beta]$,
$[x=0]$ have different meanings in Figs.\ref{fig:twoSHAs}(a),(b).
In the former, they define the guard set conditions. As confirmed
by Fig.\ref{fig:twoSHAs}(a), the guard set conditions (e.g. ,
$[\alpha>\beta]$ at $q=0$) cannot differentiate between $(i)$ a
smooth transition from the invariant $[\alpha\leq\beta]$ to
$[\alpha>\beta]$, hence, a transition to $q=1$ and $(ii)$ a jump
in $\alpha(t)$ that causes $[\alpha(t^{+})>\beta]$ to become true
without satisfying $\alpha(t)=\beta$. Recall that the former
depends on $\theta$, whereas the latter does not. On the other
hand, the set conditions on the transitions in
Fig.\ref{fig:twoSHAs}(b) have a different meaning: they identify a
condition caused by an event occurring at the same time but in a
previous transition. Thus, $[\alpha>\beta]$ is clearly associated
with a jump in $\alpha$ in the previous transition and is
independent of $\theta$. If we are to control $\theta$ to affect
the systems's performance, it is obviously crucial to identify
transitions that depend on it as opposed to ones that do not.

Another byproduct of using the guard functions is that unlike the
conventional IPA approach which normally categorizes the events
into \emph{exogenous, endogenous,} and \emph{induced} classes and
derives different equations to capture their dependence on
parameter vector $\boldsymbol{\theta}$, guard functions enable us
to treat all events uniformly. Moreover, inclusion of timers in
the states eliminates the need for a spontaneous event and the
ambiguous notions of \textquotedblleft enabled\textquotedblright\
events and \textquotedblleft waiting\textquotedblright\ in a guard
set. To show that the framework above encompasses the event
classification in \citep{CassWarPanYao09} for a SFM, we give a
simple definition of each event class in terms of guard functions:

\noindent$\bullet$ \emph{Exogenous Events}: In
\citep{CassWarPanYao09}, an event $E_{i}\in\mathcal{E}$ is defined
as exogenous if it causes a transition at time $\tau_{k}$
independent of $\boldsymbol{\theta}$ and satisfies the gradient
condition $\frac{d\tau_{k}}{d\boldsymbol{\theta}}=\mathbf{0}$. In
our case, we define $E_{i}$ occurring at $t=\tau_{k}$ as exogenous
if the associated guard function
$g_{i}(\tau_{k},\mathrm{\mathbf{x}},\mathrm{\mathbf{u}},\btheta)=g_{i}(\tau_{k},\mathrm{\mathbf{x}},\mathrm{\mathbf{u}})$
is independent of
$\boldsymbol{\theta}$ \ so that $\frac{dg_{i}(\tau_{k},\bx,\bu,\btheta)}{d\boldsymbol{\theta}%
}=\mathbf{0}$. In our framework, an exogenous event is associated with a timer
state variable whose guard function is independent of $\boldsymbol{\theta}$.

\noindent$\bullet$ \emph{Endogenous Events}: In contrast to an exogenous
event, an endogenous one is such that $\frac{dg_{i}(\tau_{k},\bx,\bu,\btheta)}%
{d\boldsymbol{\theta}}\neq\mathbf{0}$. This includes cases where $g_{i}%
(\tau_{k},\mathrm{\mathbf{x}}(\tau_{k},\boldsymbol{\theta}),\mathrm{\mathbf{u}%
}(\tau_{k},\boldsymbol{\theta}),\boldsymbol{\theta}\big)=y_{i}(\tau
_{k},\boldsymbol{\theta})=0$, for some timer state $y_{i}$.

\noindent$\bullet$ \emph{Induced Events}: In
\citep{CassWarPanYao09}, an event $E_{i}\in\mathcal{E}$ occurring
at $t=\tau_{k}$ is called \emph{induced} if it is caused by the
occurrence of another event (the \emph{triggering} event) at
$\tau_{m}<\tau_{k}$ ($m<k$). In our case, an event is induced if
there exists a state variable $x$ associated with
$E_{j}\in\mathcal{E}$ for which the following conditions are met:
\begin{subequations}
\begin{align}
\tau_{m}  &  =\max\{t<\tau_{k}:x(t,\boldsymbol{\theta})\neq0,\
\forall
t\in(\tau_{m},\tau_{k})\}.\\
\tau_{k}  &  =\min\{t>\tau_{m}:x(t,\boldsymbol{\theta})=0\}.
\end{align}
In this case, the active guard function at $\tau_{k}$ is
\end{subequations}
\begin{equation}
g_{j}(\tau_{k},\mathrm{\mathbf{x}}(\tau_{k},\boldsymbol{\theta}%
),\mathrm{\mathbf{u}}(\tau_{k},\boldsymbol{\theta}),\boldsymbol{\theta}%
)=x(\tau_{k},\boldsymbol{\theta})=0. \label{inducedg}%
\end{equation}
Generally, the initial value of state
$x(\tau_{m},\boldsymbol{\theta})$ is determined by a reset
function associated with transition $m$. After the reset, the
dynamics $\dot{x}=f(t,\boldsymbol{\theta})$ can be arbitrary until
$x(t,\boldsymbol{\theta})=0$ is satisfied. In this sense, the
timer events are simple cases of induced events with trivial
dynamics.

As already mentioned, it is possible to have simultaneous
transitions. This is a necessary condition to have chattering in
the SHA sample path which is mostly undesirable. To ensure a
bounded number of transitions in the interval $[0,T]$, let us
introduce the following assumption:

\begin{assumption}
With probability 1, the number of simultaneous transitions is
finite. \label{assump:incompatibletransitions}
\end{assumption}

Since the number of transitions occurring at different times is
finite over a finite interval $[0,T]$, Assumption
\ref{assump:incompatibletransitions} ensures that the total number
of events observed on the sample path is finite. Depending on the
system in question Assumption \ref{assump:incompatibletransitions}
translates into conditions on the states, parameters and inputs
(see Assumption \ref{assump:notequal} for the case example in
Section \ref{sec:singlenodeexample}). We do not give the
conditions under which Assumption
\ref{assump:incompatibletransitions} is valid in a general SHS
setting. More on this can be found in \citep{Sim00} and standard
references on control of hybrid systems.

\subsection{The optimization problem}

Observing the SHS just described over the interval $[0,T]$, we
seek to solve the following optimization problem:
\begin{align}
\mathbf{(P)}\quad\boldsymbol{\theta}^{\ast}  &  =\mathrm{argmin}%
_{\boldsymbol{\theta}}\ J(T,\boldsymbol{\theta})=\text{E}_{\omega
}\big[L(T,\boldsymbol{\theta},\omega)\big],\nonumber\\
\mbox{sub}  &  \mbox{ject to }\nonumber\\
&  \dot{\mathrm{\mathbf{x}}}(t,\boldsymbol{\theta},\omega)=\mathrm{\mathbf{f}%
}(q(t),t,\mathrm{\mathbf{x}}(t,\btheta,\omega),\mathrm{\mathbf{u}}(t,\btheta,\omega),\boldsymbol{\theta}%
),\nonumber\\
&  \mathrm{\mathbf{x}}\in Inv(q),\mathrm{\mathbf{u}}\in\mathcal{U}%
(\boldsymbol{\theta})\nonumber\\
&  \mathrm{\mathbf{x}}(\tau_{k}^{+},\boldsymbol{\theta},\omega
)=\mathrm{\mathbf{r}}(q_{k-1},q_{k},\mathrm{\mathbf{x}}(\tau_{k}%
,\boldsymbol{\theta},\omega),\mathrm{\mathbf{u}}(\tau_{k},\boldsymbol{\theta
},\omega),\boldsymbol{\theta})\nonumber\\
&  q(t)\in\{1,\ldots,N_{q}\},\ k=1,\ldots,K(\omega)\nonumber\\
&  \mathrm{\mathbf{x}}(0,\boldsymbol{\theta},\omega)=\mathrm{\mathbf{x}}%
_{0},\nonumber
\end{align}
where $q(t)=q_{k}$ when $t\in\lbrack\tau_{k},\tau_{k+1})$ and
$L(T,\boldsymbol{\theta},\omega)$ is a sample function generally defined as
\begin{equation}
L(T,\boldsymbol{\theta},\omega)=\int_{0}^{T}\ell(q(t),t,\mathrm{\mathbf{x}%
}(t,\btheta,\omega),\mathrm{\mathbf{u}}(t,\btheta,\omega),\boldsymbol{\theta})dt \label{samplecost}%
\end{equation}
for some given function $\ell(\cdot)$. Notice that although it is
possible to treat time $t$ as a continuous state variable, we make
the dependence of various function on $t$ explicit and do not
include it in $\bx$ as it makes our analysis easier to follow.

We solve problem (\textbf{P}) using IPA. The objective of IPA is
to specify how the changes in $\boldsymbol{\theta}$ affect
$\mathrm{\mathbf{x}}(t,\boldsymbol{\theta},\omega)$ and
ultimately, to calculate
$\frac{dL(T,\boldsymbol{\theta},\omega)}{d\boldsymbol{\theta}}$.
This is done by finding the gradients of state
$\mathrm{\mathbf{x}}$ and event times
$\tau_{k}(\boldsymbol{\theta}),$,$k=1,\ldots,K$, with respect to
$\boldsymbol{\theta}$. It has been shown that under mild
smoothness conditions the result is an unbiased estimate of the
objective function gradient
$\frac{dJ(T,\boldsymbol{\theta})}{d\boldsymbol{\theta}}$
\citep{CassWarPanYao09}. Thus, coupling the sensitivity estimates
with a gradient-based optimization algorithm can optimize the
system performance.

\section{Unified IPA Approach}

\subsection{Matrix Notation}

Let $v(t,\boldsymbol{\theta})$ be a scalar function which is differentiable
with respect to $\boldsymbol{\theta}$. We define the gradient vector with
respect to $\boldsymbol{\theta}$ as $v^{\prime}(t,\boldsymbol{\theta}%
)=\frac{\partial v(t,\boldsymbol{\theta})}{\partial\boldsymbol{\theta}%
}=\big(\frac{\partial v(t,\boldsymbol{\theta})}{\partial\theta_{1}}%
,\ldots,\frac{\partial v(t,\boldsymbol{\theta})}{\partial\theta_{N_{\theta}}%
}\big)$. Moreover, we denote the full and partial Jacobian of a vector
$\mathrm{\mathbf{v}}(t,\boldsymbol{\theta})\in\mathbb{R}^{M}$ with respect to
$\boldsymbol{\theta}$ by
\begin{align}
\frac{d\mathrm{\mathbf{v}}(t,\boldsymbol{\theta})}{d\boldsymbol{\theta}}  &
=\left[  \frac{dv_{i}(t,\boldsymbol{\theta})}{d\theta_{j}}\right]
\in\mathbb{R}^{M\times N_{\theta}},\label{fullJacobian}\\
\mathbf{\mathbf{v}^{\prime}}(t,\boldsymbol{\theta})  &  =\left[
\frac{\partial
v_{i}(t,\boldsymbol{\theta})}{\partial\theta_{j}}\right]
\in\mathbb{R}^{M\times N_{\theta}} \label{partialJacobian}%
\end{align}
where $v_{i}(t,\boldsymbol{\theta})$ ($i<M$) is the $i$-th entry of
$\mathrm{\mathbf{v}}(t,\boldsymbol{\theta})$. With a slight abuse of notation
we use $\frac{d\mathrm{\mathbf{v}}(t,\boldsymbol{\theta})}{d\mathrm{\mathbf{x}%
}}=\left[  \frac{dv_{i}(t,\boldsymbol{\theta})}{dx_{j}}\right]  $,
$\frac{\partial\mathrm{\mathbf{v}}(t,\boldsymbol{\theta})}{\partial
\mathrm{\mathbf{x}}}=\left[  \frac{\partial v_{i}(t,\boldsymbol{\theta}%
)}{\partial x_{j}}\right]  \in\mathbb{R}^{M\times N_{x}}$ and $\frac
{d\mathrm{\mathbf{v}}(t,\boldsymbol{\theta})}{d\mathrm{\mathbf{u}}}=\left[
\frac{dv_{i}(t,\boldsymbol{\theta})}{du_{k}}\right]  $, $\frac{\partial
\mathrm{\mathbf{v}}(t,\boldsymbol{\theta})}{\partial\mathrm{\mathbf{u}}%
}=\left[  \frac{\partial v_{i}(t,\boldsymbol{\theta})}{\partial u_{k}}\right]
\in\mathbb{R}^{M\times N_{u}}$ as the full and partial Jacobians of
$\mathrm{\mathbf{v}}(t,\boldsymbol{\theta})$ with respect to
$\mathrm{\mathbf{x}}$ and $\mathrm{\mathbf{u}}$.

For the event times $\tau_{k}(\boldsymbol{\theta})$, $k=1,\ldots,K$, the
gradient with respect to $\boldsymbol{\theta}$ is defined as
\[
\tau_{k}^{\prime}=(\tau_{k,1}^{\prime},\ldots,\tau_{k,N_{\theta}}^{\prime})
\]
where $\tau_{k,j}^{\prime}=\frac{\partial\tau_{k}}{\partial\theta_{j}}$. We
let $\tau_{0,j}^{\prime}=\tau_{K+1,j}^{\prime}=0$ for all $j$ since the start
and end of the sample path are fixed values. Finally, we define
$\boldsymbol{\tau}^{\prime}$ as a $N_{e}\times N_{\theta}$ matrix such that
its $i$th row is associated with event $E_{i}$ and its $j$th column is
associated with the variable with respect to which the differentiation is done.

In what follows, we derive a unified set of equations which give the
event-time and state derivatives with respect to $\boldsymbol{\theta}$ and are
in concord with the results of IPA presented in \citep{CassWarPanYao09}. All
calculations can be done in two generic steps, regardless of the type of event
observed, i.e., we do not need to differentiate between exogenous, endogenous,
and induced events.

\section{Infinitesimal Perturbation Analysis}

Below, we drop $T$ from $L(T,\btheta,\omega)$ and $\omega$ from
the arguments of other functions to simplify the notation.
However, we still write $L(\boldsymbol{\theta},\omega)$ to stress
that we carry out the analysis on the sample path of system $G$
denoted by $\omega$. We write (\ref{samplecost}) as
\[
L(\boldsymbol{\theta},\omega)=\sum_{k=0}^{K}\int_{\tau_{k}}^{\tau_{k+1}}%
\ell(q_{k},t,\mathrm{\mathbf{x}},\mathrm{\mathbf{u}},\boldsymbol{\theta
})dt,
\]
Recalling $\tau_{0}^{\prime}=\tau_{K+1}^{\prime}=\mathbf{0}$, we calculate the
gradient of the sample cost with respect to $\boldsymbol{\theta}$ as
\begin{align}
\frac{dL(\boldsymbol{\theta},\omega)}{d\boldsymbol{\theta}}  &  =\sum
_{k=1}^{K}[\ell(q_{k-1},\tau_{k}^{-},\mathrm{\mathbf{x}},\mathrm{\mathbf{u}%
},\boldsymbol{\theta})-\ell(q_{k},\tau_{k}^{+},\mathrm{\mathbf{x}%
},\mathrm{\mathbf{u}},\boldsymbol{\theta})]\tau_{k}^{\prime}\nonumber\\
&  +\sum_{k=0}^{K}\int_{\tau_{k}}^{\tau_{k+1}}\frac{d\ell(q_{k}%
,t,\mathrm{\mathbf{x}},\mathrm{\mathbf{u}},\boldsymbol{\theta})}%
{d\boldsymbol{\theta}}dt \label{dLdtheta}%
\end{align}
where
\begin{align}
\frac{d\ell(q_{k},t,\mathrm{\mathbf{x}},\mathrm{\mathbf{u}},\boldsymbol{\theta
})}{d\boldsymbol{\theta}}=\frac{\partial\ell}{\partial
\mathrm{\mathbf{x}}}\frac{d\mathrm{\mathbf{x}}(t,\boldsymbol{\theta})}{d\boldsymbol{\theta}}  &  +\frac{\partial\ell}{\partial\mathrm{\mathbf{u}}}\mathrm{\mathbf{u}}^{\prime}%
(t,\boldsymbol{\theta})\nonumber\\
&  +\ell^{\prime}(q_{k},t,\mathrm{\mathbf{x}},\mathrm{\mathbf{u}%
},\boldsymbol{\theta})\nonumber
\end{align}
Thus, as mentioned before, in order to determine the sample cost gradient with
respect to $\boldsymbol{\theta}$, one needs to find the event time and state
derivatives with respect to it.

\subsection{Event-time Derivatives}

By Theorem \ref{thm:existence}, for any $k=1,\ldots,K$, transition
$k$ is directly or indirectly associated with an event
$e\in\tilde{\mathcal{E}}$ with a guard function $g(\cdot)$ such
that $g(\tau_{k}^{-},\tilde\bx,\bu,\btheta)=0$. Let us define the
guard vector
\begin{equation}
{\mathbf{g}}(t,\mathrm{\mathbf{x}},\mathrm{\mathbf{u}},\boldsymbol{\theta
})=\big(g_{1}(t,\mathrm{\mathbf{x}},\mathrm{\mathbf{u}},\boldsymbol{\theta
}),\ldots,g_{N_{e}}(t,\mathrm{\mathbf{x}},\mathrm{\mathbf{u}}%
,\boldsymbol{\theta})\big). \label{guardvector}%
\end{equation}
and a unit \emph{firing vector} $\mathrm{\mathbf{e}}_{i}$, $i=1,\ldots,N_{e}$
as
\[
\mathrm{\mathbf{e}}_{i}=\big(0,\ldots,1,\ldots,0\big)\in\mathbb{R}^{N_{e}}%
\]
where only the $i$-th element is 1 and the rest are 0.

Let ${\mathrm{\mathbf{G}}}(t,\bx,\bu,\boldsymbol{\theta})$ be a
diagonal matrix function where $\mathrm{\mathbf{G}}_{i,i}=g_{i}$,
$i=1,\ldots,N_{e}$ and denote its time derivative at $t=\tau$ by
$\dot{\mathrm{\mathbf{G}}}(\tau ,\bx,\bu,\boldsymbol{\theta})$. We
can obtain $\boldsymbol{\tau}^{\prime}$ by
differentiating ${\mathbf{g}}(\tau,\mathrm{\mathbf{x}},\mathrm{\mathbf{u}%
},\boldsymbol{\theta})=\mathbf{0}$ with respect to $\boldsymbol{\theta}$. This
gives
\begin{equation}
\boldsymbol{\tau}^{\prime}=-\dot{\mathrm{\mathbf{G}}}(\tau^{-}%
,\bx,\bu,\boldsymbol{\theta})^{-1}\frac{d{\mathbf{g}}(\tau^{-},\bx,\bu,\boldsymbol{\theta}%
)}{d\boldsymbol{\theta}} \label{gentkprime}%
\end{equation}
where
\begin{align}
\frac{d{\mathbf{g}}(\tau^{-},\bx,\bu,\boldsymbol{\theta})}{d\boldsymbol{\theta}}
&=\frac{\partial{\mathbf{g}}}{\partial\mathrm{\mathbf{x}}}\frac{d\bx(\tau^-,\btheta)}{d\btheta}\nonumber\\&+\frac{\partial{\mathbf{g}}}%
{\partial\mathrm{\mathbf{u}}}\mathrm{\mathbf{u}}^{\prime}(\tau^{-}%
,\boldsymbol{\theta})+{\mathbf{g}}^{\prime}(\tau^-,\mathrm{\mathbf{x}}%
,\mathrm{\mathbf{u}},\boldsymbol{\theta}). \label{dbgdtheta}%
\end{align}
It is easily verified that the simple equation (\ref{gentkprime})
is in line with what has been reported in prior work on IPA for
SHS, e.g., in \citep{CassWarPanYao09}. The following assumption is
introduced so that $\tau_{k}^{\prime}$ exists:

\begin{assumption}
With probability 1, if $\tau_{k}$ is the occurrence time of
$E_{i}$, we have
$\dot{g}_{i}(\tau_{k},\bx,\bu,\boldsymbol{\theta})\neq0$.
\label{assump:donttouch}
\end{assumption}

Note that the case of a contact point where $\dot{g}_{i}(\tau_{k}%
,\boldsymbol{\theta})$ does not exist has already been excluded by
Assumption \ref{assump:nosimul}, hence,
$\dot{g}_{i}(\tau_{k},\bx,\bu,\boldsymbol{\theta})$ is always
well-defined. Also, observe that only one row in
(\ref{gentkprime}) is evaluated at each transition. In fact,
$\boldsymbol{\tau}^{\prime}$ is a generic matrix function such
that if transition $k$ is associated with event $E_{i}$, only its
$i$th row is evaluated.

\subsection{State Derivatives}

Here, we determine how the state derivatives evolve both at
transition times and in between them, i.e., within a mode $q\in
Q$.

\subsubsection{Derivative update at transition times:}

At each transition we consider two cases:

$(a)$ \emph{No Reset}: In this case, we have $x_{j}(\tau_{k}^{-}%
,\boldsymbol{\theta})=x_{j}(\tau_{k}^{+},\boldsymbol{\theta})$ for all
$j\in\{1,\ldots,N_{x}\}$. Assume that for every $q\in Q$, $\dot{x}%
_{j}(t,\boldsymbol{\theta})=f_{j}(q,t,\mathrm{\mathbf{x}},\mathrm{\mathbf{u}%
},\boldsymbol{\theta})$. Then, we use the following equation, derived in
\citep{CassWarPanYao09}, to update the state derivatives:
\begin{align}
x_{j}^{\prime}&(\tau_{k}^{+},\boldsymbol{\theta})=x_{j}^{\prime}(\tau_{k}%
^{-},\boldsymbol{\theta})\nonumber\\&+[f_{j}(q_{k-1},\tau_{k}^{-},\bx,\bu,\boldsymbol{\theta
})-f_{j}(q_{k},\tau_{k}^{+},\bx,\bu,\boldsymbol{\theta})]\tau_{k}^{\prime}.
\label{update1}%
\end{align}

$(b)$ \emph{Reset}: In this case, there exists $j\in\{1,\ldots,N_{x}\}$ such
that $x_{j}(\tau_{k}^{-},\boldsymbol{\theta})\neq x_{j}(\tau_{k}%
^{+},\boldsymbol{\theta})$. Let $e\in\mathcal{E}$ be the direct cause of
transition $(q_{k-1},q_{k})$ at $\tau_{k}$ (i.e., $e$ appears on the arc
connecting $q_{k-1}$ and $q_{k}$ in the automaton). We then define a reset
condition $x_{j}(\tau_{k}^{+},\boldsymbol{\theta})=r_{j}(q_{k-1}%
,q_{k},\mathrm{\mathbf{x}},\mathrm{\mathbf{u}},\boldsymbol{\theta},e)$
where $r_{j}(\cdot)$ is the reset function of $x_{j}$. Thus, we
get
\begin{equation}
x_{j}^{\prime}(\tau_{k}^{+},\boldsymbol{\theta})=\frac{dr_{j}^{\prime}%
(q_{k-1},q_{k},\mathrm{\mathbf{x}},\mathrm{\mathbf{u}},\boldsymbol{\theta}%
,e)}{d\boldsymbol{\theta}} \label{update2}%
\end{equation}
where
\begin{align}
\frac{dr_{j}^{\prime}(q_{k-1},q_{k},\mathrm{\mathbf{x}},\mathrm{\mathbf{u}%
},\boldsymbol{\theta},e)}{d\boldsymbol{\theta}}=  &  \frac{\partial r_{j}%
}{\partial\mathrm{\mathbf{x}}}\mathrm{\mathbf{x}}^{\prime}(\tau_{k}%
,\boldsymbol{\theta})+\frac{\partial r_{j}}{\partial\mathrm{\mathbf{u}}%
}\mathrm{\mathbf{u}}^{\prime}(\tau_{k},\boldsymbol{\theta})\nonumber\\
&  +r_{j}^{\prime}(q_{k-1},q_{k},\mathrm{\mathbf{x}},\mathrm{\mathbf{u}%
},\boldsymbol{\theta},e). \label{dresetdtheta}%
\end{align}
For other transitions which are indirectly caused by an event, we simply
define $\mathrm{\mathbf{r}}(q_{k-1},q_{k},\mathrm{\mathbf{x}}%
,\mathrm{\mathbf{u}},\boldsymbol{\theta})=\mathrm{\mathbf{x}}(\tau_{k},\btheta)$
as no reset is possible on them.

To put everything in matrix form, let us first define, for every
$i=1,\ldots,N_{e}$, the index set
\[
\Phi_{i}=\{(m,n)\in Q\times Q:\exists\mathrm{\mathbf{x}},\mathrm{\mathbf{u}%
}\in\mathcal{X},\mathcal{U}\ \mbox{s.t.}\phi(m,\mathrm{\mathbf{x}%
},\mathrm{\mathbf{u}},E_{i})=n\}
\]
containing all transitions directly associated with $E_{i}$. Also, for each
$i$, let us define the reset mapping
\[
\mathrm{\mathbf{r}}_{i}(m,n)=\left\{
\begin{array}
[c]{ll}%
\mathrm{\mathbf{r}}(m,n,\mathrm{\mathbf{x}},\mathrm{\mathbf{u}}%
,\boldsymbol{\theta},E_{i}) & \mathrm{\ if\ }(m,n)\in\Phi_{i}\\
\mathrm{\mathbf{x}} & \mathrm{\ otherwise}%
\end{array}
\right.
\]
and the diagonal matrix $\mathrm{\mathbf{C}}(m,n)\in\mathbb{R}^{N_{x}\times
N_{x}}$ with its $j$th diagonal entry $c_{jj}(m,n)=1$ if $x_{j}$ is not reset
by transition $(m,n)$ and 0, otherwise. We also define $\bar{\mathrm{\mathbf{C}%
}}(m,n)=\mathrm{\mathbf{I}}_{N_{x}\times N_{x}}-\mathrm{\mathbf{C}}(m,n)$
where $\mathrm{\mathbf{I}}_{N_{x}\times N_{x}}$ is the identity matrix with
the specified dimensions. Moreover, let us define the reset map matrix as
\begin{equation}
\mathrm{\mathbf{R}}(\tau_{k})=%
\begin{bmatrix}
\mathrm{\mathbf{r}}_{1}(q_{k-1},q_{k})\\
\vdots\\
\mathrm{\mathbf{r}}_{N_{e}}(q_{k-1},q_{k})
\end{bmatrix}
\in\mathbb{R}^{N_{e}\times N_{x}} \label{resetmatrix}%
\end{equation}
where $x_{j}(\tau_{k}^{+})=r_{i,j}(q_{k-1},q_{k})$ when the $k$th transition
is due to $E_{i}$ and $x_{j}(\tau_{k}^{+})=x_{j}(\tau_{k}^{-})$, otherwise.
Thus, if $\mathrm{\mathbf{x}}$ remains continuous at its occurrence time we
get $\mathrm{\mathbf{r}}_{i}(q_{k-1},q_{k})=\mathrm{\mathbf{x}}(\tau
_{k},\boldsymbol{\theta})$. Finally, we define the shorthand notation
$\Delta\mathrm{\mathbf{f}}(q_{k-1},q_{k},\boldsymbol{\theta}%
)\equiv\mathrm{\mathbf{f}}(q_{k-1},\tau_{k}^{-},\mathrm{\mathbf{x}}%
,\mathrm{\mathbf{u}},\boldsymbol{\theta})-\mathrm{\mathbf{f}}(q_{k},\tau
_{k}^{+},\mathrm{\mathbf{x}},\mathrm{\mathbf{u}},\boldsymbol{\theta})$ for the
jump in the dynamics at the $k$th transition. Using these definitions we can
combine part $(a)$ and $(b)$ above and write
\begin{align}
\mathrm{\mathbf{x}}^{\prime}(\tau_{k}^{+})=\mathrm{\mathbf{C}}(q_{k-1},q_{k})
&  \big[\mathrm{\mathbf{x}}^{\prime}(\tau_{k}^{-})+\Delta\mathrm{\mathbf{f}%
}(q_{k-1},q_{k},\boldsymbol{\theta})^{\mathsf{T}}\tau_{k}^{\prime
}\big]\nonumber\\
&  +\left(  \mathrm{\mathbf{e}}_{i}\mathrm{\mathbf{R}}(\tau_{k})\bar
{\mathrm{\mathbf{C}}}(q_{k-1},q_{k})\right)  ^{\prime}. \label{discupdate}%
\end{align}

\subsubsection{Derivative update between transitions:}

Assuming the mode is $q_{k}$, we only need to perform the following operation
on interval $[\tau_{k},\tau_{k+1})$:%
\begin{equation}
\mathrm{\mathbf{x}}^{\prime}(t,\boldsymbol{\theta})=\mathrm{\mathbf{x}%
}^{\prime}(\tau_{k}^{+},\boldsymbol{\theta})+\int_{\tau_{k}}^{t}%
\frac{d\mathrm{\mathbf{f}}(q_{k},\tau,\mathrm{\mathbf{x}},\mathrm{\mathbf{u}%
},\boldsymbol{\theta})}{d\boldsymbol{\theta}}d\tau \label{contupdate}%
\end{equation}
where $\frac{d\mathrm{\mathbf{f}}(q_{k},\tau,\mathrm{\mathbf{x}},\mathrm{\mathbf{u}%
},\boldsymbol{\theta})}{d\boldsymbol{\theta}}$ is a $N_{x}\times
N_{\theta}$ Jacobian matrix of the state dynamics defined on
$[\tau_{k},\tau_{k+1})$ as
\begin{align}
\frac{d\mathrm{\mathbf{f}}(q_{k},\tau,\mathrm{\mathbf{x}},\mathrm{\mathbf{u}%
},\boldsymbol{\theta})}{d\boldsymbol{\theta}}=\frac{\partial\mathrm{\mathbf{f}%
}}{\partial\mathrm{\mathbf{x}}}  &  \frac{d\bx(\tau,\btheta)}{d\btheta}+\frac{\partial\mathrm{\mathbf{f}}}%
{\partial\mathrm{\mathbf{u}}}\mathrm{\mathbf{u}}^{\prime}(\tau
,\boldsymbol{\theta})\nonumber\\
&  +\mathrm{\mathbf{f}}^{\prime}(q_{k},\tau,\mathrm{\mathbf{x}}%
,\mathrm{\mathbf{u}},\boldsymbol{\theta}). \label{dfdtheta}%
\end{align}


To sum up, the basis for IPA on a system modeled as $G$ is the
pre-calculation of the quantities $\frac{\partial\ell}{\partial
\mathrm{\mathbf{x}}},\frac{\partial\ell}{\partial\mathrm{\mathbf{u}}},$
$\frac{\partial\mathrm{\mathbf{f}}}{\partial\mathrm{\mathbf{x}}}$,
$\frac{\partial\mathrm{\mathbf{f}}}{\partial\mathrm{\mathbf{u}}}$
for all
$q\in Q$, $\frac{\partial\mathrm{\mathbf{r}}_{i}}{\partial\mathrm{\mathbf{x}}%
}$,
$\frac{\partial\mathrm{\mathbf{r}}_{i}}{\partial\mathrm{\mathbf{u}}}$
and $\frac{\partial{\mathbf{g}}}{\partial\boldsymbol{\theta}}$,
which are then used in (\ref{gentkprime}), (\ref{discupdate}) and
(\ref{contupdate}) to update IPA derivatives. Finally, the results
are applied to (\ref{dLdtheta}). We will next apply this method to
a specific problem of interest \citep{cassbook} in the following
section.

\section{IPA for a single-node SFM}\label{sec:singlenodeexample}

\begin{figure}[th]
\centering
\includegraphics[scale =.28]{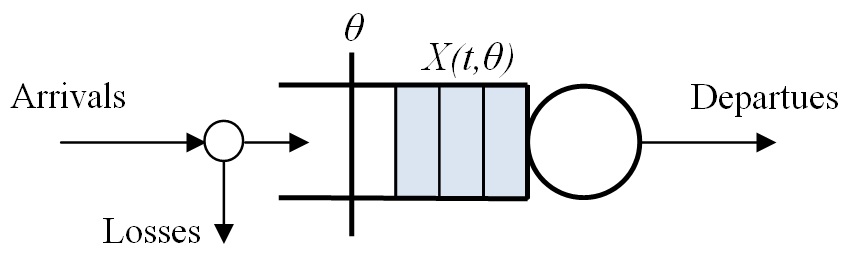}\caption{{\small The
DES model of the single node system.}}%
\label{fig:singleclassqueue}\vspace{-2mm}
\end{figure}

In what follows, we apply the method described above to a simple single-class
single-node system shown in Fig.\ref{fig:singleclassqueue}.  We use a simplified notation here for space limitations. However, the dependence of functions on their arguments should be clear from the analysis in the previous section.

The system consists of a queue whose content level $X(t,\theta)$
is subject to stochastic arrival and service time processes. The
queue capacity is limited to a quantity $\theta$ treated as the
control parameter. Every arrival seeing a full queue is lost and
incurs a penalty. Considering this system over a finite interval
$[0,T]$, we want to find the best $\theta$ to trade off between
the average workload and average loss defined as
\begin{align}
\text{E}_{\omega}[Q_{DES}(T,\theta,\omega)]  &  =\frac{1}{T}\text{E}_{\omega
}\left[  \int_{0}^{T}X(t,\theta,\omega)dt\right], \label{DESQ}\\
\text{E}[L_{DES}(T,\theta,\omega)]  &  =\frac{1}{T}\text{E}_{\omega}%
N_{loss}(T,\theta,\omega) \label{DESL}%
\end{align}
where $N_{loss}(\cdot)$ is the number of losses observed in the
interval $[0,T]$. Even for a simple system like this, the analysis
can become prohibitive when the stochastic processes considered
are arbitrary. Use of SFMs has proven to be very helpful in the
analysis and optimization of queuing systems such as this one
\citep{cassbook},\citep{CassWarMelSunPan02} where applying IPA has
resulted in very simple derivative estimates of the loss and
workload objectives with respect to $\theta$. In the SFM, the
arrivals and departures are abstracted into non-negative
stochastic inflow rate $\{\alpha(t)\}$ and maximal service rate
$\{\beta(t)\}$ processes which are independent of $\theta$. These
rates continuously evolve according to the differential equations
$\dot{\alpha }=f_{\alpha}(t)$ and $\dot{\beta}=f_{\beta}(t)$ where
$f_{\alpha}$ and $f_{\beta}$ are arbitrary continuous functions.
It is important to observe that the precise nature of these
functions turns out to be irrelevant as far as the resulting IPA
estimator is concerned: the IPA estimators are independent of
$f_{\alpha}$ and $f_{\beta}$. This is an important robustness
property of IPA estimators which holds under certain conditions
\citep{YaoCass11}. We allow for discrete jumps in both processes
$\{\alpha(t)\}$ and $\{\beta(t)\}$ and use timer states
$y_{\alpha},y_{\beta }\geq0$ to capture them. The fluid discharge
rate $d(t,\theta)$ is defined as $d(t,\theta)=\beta(t)$ when
$x(t,\theta)>0$ and $d(t,\theta)=\alpha(t)$ otherwise. For this
example, Assumption \ref{assump:incompatibletransitions} manifests
itself as follows:

\begin{assumption}
With probability 1, condition $\alpha(t)=\beta(t)$ cannot be valid
on a non-empty interval containing $t$. \label{assump:notequal}
\end{assumption}

The buffer content process evolves according to the differential equation
$\dot{x}(t,\theta)=\alpha(t)-d(t,\theta)$ so we can write
\begin{equation}\label{dx}
\dot{x}(t,\theta)=f_x(t,\theta)=\left\{
\begin{array}
[c]{ll}%
0 & \mathrm{\ if\ }x(t,\theta)=0\mathrm{\ or\ }\theta\\
\alpha(t)-\beta(t) & \mathrm{\ otherwise}%
\end{array}
\right.
\end{equation}
When $x(t,\theta)$ reaches the buffer capacity level $\theta$, a portion of
the incoming flow is rejected with rate $\alpha(t)-\beta(t)>0$. Obviously,
when $x(t,\theta)<\theta$ no loss occurs. Hence, we define, for every
$t\in\lbrack0,T]$, the \emph{loss rate} as
\begin{equation}
\ell(t,\theta)=\left\{
\begin{array}
[c]{ll}%
0 & \mathrm{\ if\ }x(t,\theta)<\theta\\
\alpha(t)-\beta(t) & \mathrm{\ otherwise}%
\end{array}
\right.  \label{lossrate}%
\end{equation}
\begin{figure}[th]
\centering
\includegraphics[scale =.175]{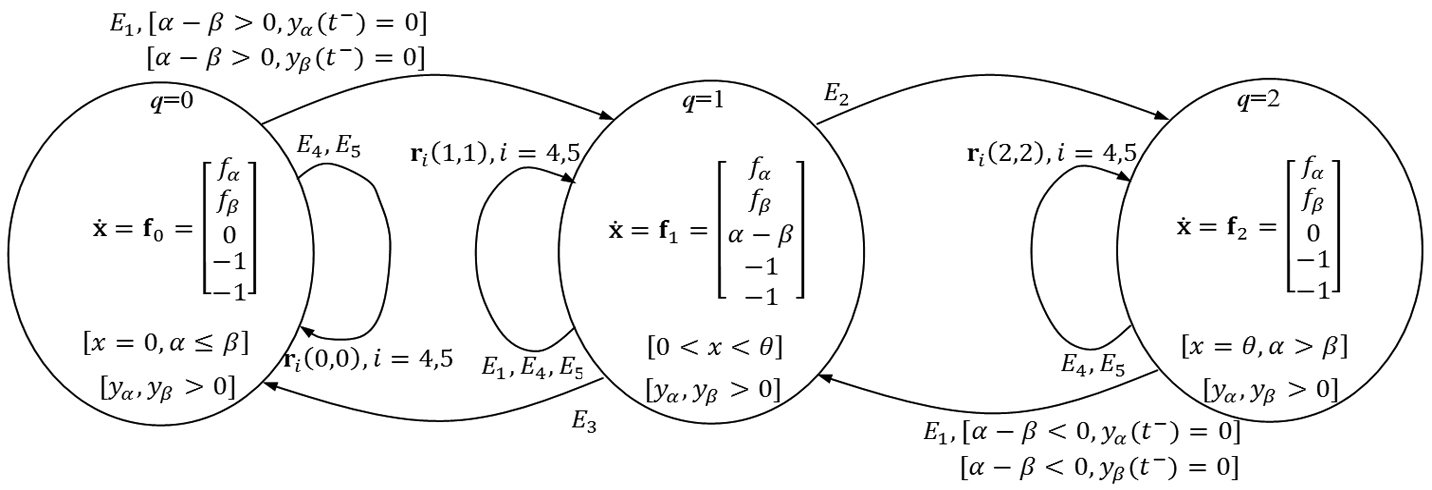}\caption{{\small Stochastic
Flow Model (SFM) for the single-class SFM.}}%
\label{fig:singleclassSFM}\vspace{-2mm}
\end{figure}The SHA of this system is shown in Fig. \ref{fig:singleclassSFM}.
We define the system state vector $\mathrm{\mathbf{x}}=(\alpha,\beta
,x,y_{\alpha},y_{\beta})$ and the input $\mathrm{\mathbf{u}}=(\{A_{j}%
\},\{V_{j}^{\alpha}\},\{B_{k}\},\{V_{k}^{\beta}\})$ whose elements are
sequences of random variables from the jump distributions associated with
states $\alpha$, $y_{\alpha}$, $\beta$, and $y_{\beta}$, respectively.
Although it follows from the definitions that the states $\alpha
,\ \beta,\ y_{\alpha},\ y_{\beta}$ are independent of $\theta$ yielding
$\alpha^{\prime}(t)=\beta^{\prime}(t)=y_{\alpha}^{\prime}(t)=y_{\beta}%
^{\prime}(t)=0$ for all $t\in\lbrack0,T]$, we include them at the start of the
IPA estimation procedure to fully illustrate the matrix setting we have
developed and make use of it later. Regarding the dynamics of the states, we
have
\begin{equation}
\mathrm{\mathbf{f}}_{q}=\left\{
\begin{array}
[c]{ll}%
(f_{\alpha},f_{\beta},0,-1,-1) & \mathrm{\ if\ }q=0\\
(f_{\alpha},f_{\beta},\alpha-\beta,-1,-1) & \mathrm{\ if\ }q=1\\
(f_{\alpha},f_{\beta},0,-1,-1) & \mathrm{\ otherwise}%
\end{array}
\right. \label{fq}\end{equation} The workload and loss sample
functions are
\begin{align}
Q(T,\theta,\omega)  &  =\frac{1}{T}\int_{0}^{T}x(t,\theta,\omega)dt,\nonumber\\
L(T,\theta,\omega)  &
=\frac{1}{T}\int_{0}^{T}\ell(t,\theta,\omega )dt.\nonumber
\end{align}
Using the fact that $x(t,\theta)$ and $\ell(t)$ can only contribute to their
associated objectives when, respectively, $x>0$ ($q=1,2$) and $x=\theta$
($q=2$), we can write:
\begin{align}
Q(T,\theta,\omega)  &  =\frac{1}{T}\sum_{n=1}^{N}\int_{\xi_{n}(\theta)}%
^{\eta_{n}(\theta)}x(t,\theta,\omega)dt,\label{Q}\\
L(T,\theta,\omega)  &  =\frac{1}{T}\sum_{n=1}^{N}\sum_{m=1}^{M_{n}}\int
_{\nu_{n,m}(\theta)}^{\sigma_{n,m}(\theta)}\ell(t,\theta,\omega)dt \label{L}%
\end{align}
where $N$ is the number of supremal intervals $[\xi_{n},\eta_{n}),$
$n=1,\ldots,N$ over which $x(t,\theta,\omega)>0$ (i.e., $q=1,2$) and $M_{n}$
is the number of supremal intervals $[\nu_{n,m},\sigma_{n,m}),$ $m=1,\ldots
,M_{n}$ such that $x(t,\theta,\omega)=\theta$ $($i.e., $q=2)$. We refer to the
intervals of the first kind as \emph{Non-Empty Periods} (NEPs) and the second
kind as \emph{Full Periods} (FPs). We also drop the sample path index $\omega$
to simplify the notation.

Differentiating \eqref{Q} with respect to $\theta$ and noting $x(\xi
_{n}(\theta)^{+})=x(\eta_{n}(\theta)^{-})=0$ for any $n$ gives
\begin{align}
&  Q^{\prime}(T,\theta)=\frac{1}{T}\sum_{n=1}^{N}\bigg[
x(\eta_{n}(\theta)^{-})\eta_n^\prime(\theta)-x(\xi_{n}(\theta
)^{+})\xi_n^\prime(\theta)\nonumber\\
&  \
+\int_{\xi_{n}(\theta)}^{\eta_{n}(\theta)}x^{\prime}(t,\theta)dt\bigg]=
\frac{1}{T}\sum_{n=1}^{N}\int_{\xi_{n}(\theta)}^{\eta_{n}(\theta)}x^{\prime
}(t,\theta)dt, \label{Qprime}%
\end{align}
Next, differentiating \eqref{L} with respect to $\theta$ and noticing that by
(\ref{lossrate}), $\ell^{\prime}(t,\theta)=0$ for all $t\in\lbrack0,T)$
(eliminating the integral part), reveals
\begin{equation}
L^{\prime}(T,\theta)=\sum_{n=1}^{N}\sum_{m=1}^{M_{n}}[\sigma_{n,m}^{\prime
}\ell(\sigma_{n,m}^{-},\theta)-\nu_{n,m}^{\prime}\ell(\nu_{n,m}^{+},\theta)].
\label{Lprime}%
\end{equation}
It is now clear that to evaluate \eqref{Qprime} and \eqref{Lprime} we only
need to obtain $x^{\prime}(t,\theta)$ for all $t\in\lbrack\xi_{n},\eta_{n})$,
$n=1,\ldots,N$ and event time derivatives $\nu_{n,m}^{\prime},\sigma
_{n,m}^{\prime}$ for every $m=1,\dots,M_{n}$ where $n=1,\ldots,N$. However, as
mentioned before, we try to keep everything in the general matrix framework so
as to verify its effectiveness.

According to Fig. \ref{fig:singleclassSFM}, the event set of the SHS is given
as
\[
\mathcal{E}=\{E_{i},i=1,\ldots,5\}.
\]
with guard functions defined as follows: $E_{1}$ occurs when
$g_{1}=\alpha-\beta=0$. $E_{2}$ is the event of reaching the
buffer threshold $\theta$, so that $g_{2}=x-\theta$. $E_{3}$ is
the event ending a non-empty period, hence $g_{3}=x$. Finally,
$E_{4}$ and $E_{4}$ are associated with the timer run-offs
captured by $g_{4}=y_{\alpha}$ and $g_{5}=y_{\beta}$,
respectively. To summarize, the guard vector for the system is
\begin{equation}
{\mathbf{g}}(t,\theta)=\big(\alpha(t)-\beta(t),x(t,\theta)-\theta
,x(t,\theta),y_{\alpha}(t),y_{\beta}(t)\big). \label{guards}%
\end{equation}
The reset maps are defined as follows:
\[
\mathrm{\mathbf{r}}_{i}(m,n)=\vspace{-2mm}\left\{
\begin{array}
[c]{ll}%
\hspace{-1.5mm}(A_{j},\beta,x,V_{j}^{\alpha},y_{\beta}) & \mathrm{\ if\ }%
i=4,m=n\in\{0,1,2\}\\
\hspace{-1.5mm}(\alpha,B_{k},x,y_{\alpha},V_{k}^{\beta}) & \mathrm{\ if\ }%
i=5,m=n\in\{0,1,2\}\\
\hspace{-1.5mm}\mathrm{\mathbf{x}} & \mathrm{\ otherwise}%
\end{array}
\right.
\]
where $A_{j}$ and $B_{k}$ are, respectively, the $j$th and $k$th elements of
random sequences $\{A_{j}\}$ and $\{B_{k}\}$. Using these results in
(\ref{resetmatrix}) we get
\begin{equation}
\mathrm{\mathbf{R}}(\tau_{k})=%
\begin{bmatrix}
\mathrm{\mathbf{x}}(\tau_{k})\\
\mathrm{\mathbf{x}}(\tau_{k})\\
\mathrm{\mathbf{x}}(\tau_{k})\\
\mathrm{\mathbf{r}}_{4}(q_{k-1},q_{k})\\
\mathrm{\mathbf{r}}_{5}(q_{k-1},q_{k})
\end{bmatrix}
. \label{exampleresetmatrix}%
\end{equation}
By (\ref{discupdate}), only the reset conditions associated with discontinuous
states need to be differentiated with respect to $\theta$. Since $x$, the only
state variable which depends of $\theta$, is continuous, we conclude that the
last term in (\ref{discupdate}) is always $\mathbf{0}$ and need not be evaluated.

The IPA starts by evaluating (\ref{gentkprime}). Note that by
(\ref{fq}) $\frac
{d\mathrm{\mathbf{f}}_{q}(t)}{d\theta}=\mathbf{0}$ for all
$q\in\{0,1,2\}$. Thus, we only need to evaluate (\ref{dbgdtheta}).
Also, by definition,
$\frac{\partial\mathrm{\mathbf{u}}}{\partial\theta}=\mathbf{0}$,
so $\frac{\partial{\mathbf{g}}}{\partial\mathrm{\mathbf{u}}}$ need
not be evaluated. Then, from (\ref{guards}), we are left with
\[
\frac{\partial{\mathbf{g}}}{\partial\mathrm{\mathbf{x}}}=%
\begin{bmatrix}
1 & -1 & 0 & 0 & 0\\
0 & 0 & 1 & 0 & 0\\
0 & 0 & 1 & 0 & 0\\
0 & 0 & 0 & 1 & 0\\
0 & 0 & 0 & 0 & 1
\end{bmatrix}
,\ {\mathbf{g}}^{\prime}(t,\theta)=%
\begin{bmatrix}
0\\
-1\\
0\\
0\\
0
\end{bmatrix}
.
\]
It follows from (\ref{dbgdtheta}) that%

\[
\frac{d{\mathbf{g}}(\tau^{-},\theta)}{d\theta}=%
\begin{bmatrix}
\alpha^{\prime}(\tau^{-})-\beta^{\prime}(\tau^{-})\\
x^{\prime}(\tau^{-})-1\\
x^{\prime}(\tau^{-})\\
y_{\alpha}^{\prime}(\tau^{-})\\
y_{\beta}^{\prime}(\tau^{-})
\end{bmatrix}
=%
\begin{bmatrix}
0\\
x^{\prime}(\tau^{-})-1\\
x^{\prime}(\tau^{-})\\
0\\
0
\end{bmatrix}
\]
where we used the fact that $\alpha$, $\beta,y_{\alpha}$ and $y_{\beta}$ are
independent of $\theta$. Moreover, by (\ref{guards}), we have the time
derivative of the guard functions just before the $k$th transition as
\[
\dot{\mathrm{\mathbf{G}}}(\tau^{-},\theta)=\mathrm{diag}\big(f_{\alpha}%
(\tau^{-})-f_{\beta}(\tau^{-}),\dot{x}(\tau^{-}),\dot{x}(\tau^{-}%
),-1,-1\big)
\]
and since $E_{2}$ and $E_{3}$ are only feasible when $x>0$ ($q=1,2$), by
(\ref{dx}) we get $\dot{x}(\tau^{-})=\alpha(\tau^{-})-\beta(\tau^{-})\neq0$ in
the above expression. Combining the results in (\ref{gentkprime}) gives
\begin{equation}
\boldsymbol{\tau}^{\prime}=-\left(  0,\frac{x^{\prime}(\tau^{-})-1}%
{\alpha(\tau^{-})-\beta(\tau^{-})},\frac{x^{\prime}(\tau^{-})}{\alpha(\tau
^{-})-\beta(\tau^{-})},0,0\right)^{\mathsf{T}}.\hspace{-2mm}
\label{finaltkprime}%
\end{equation}
Next, we determine the state derivatives with respect to $\theta$.
Since $x$ is the only state dependent on $\theta$, we only apply
(\ref{discupdate}) and (\ref{contupdate}) to $x(t,\theta)$. We
also need not evaluate $\mathrm{\mathbf{C}}(q_{k-1},q_{k})$,
$k=1,\ldots,K$ as $x(t,\theta)$ is continuous
throughout $[0,T]$. Moreover, we need to determine $\Delta f_{x}(q_{k-1}%
,q_{k})$ for any feasible transition $(q_{k-1},q_{k})$. If we define
$\lambda(\tau)=\alpha(\tau)-\beta(\tau)$ and $\Delta\lambda(\tau)=\lambda
(\tau^{-})-\lambda(\tau^{+})$ we get
\[
\Delta f_{x}(q_{k-1},q_{k})=\hspace{-1mm}\left\{  \hspace{-1.5mm}%
\begin{array}
[c]{ll}%
-\lambda(\tau_{k}^{+}) & \mathrm{\ if\ }(q_{k-1},q_{k})\in\{(0,1),(2,1)\}\\
\lambda(\tau_{k}^{-}) & \mathrm{\ if\ }(q_{k-1},q_{k})\in\{(1,0),(1,2)\}\\
\Delta\lambda(\tau_{k}) & \mathrm{\ if\ }(q_{k-1},q_{k})=(1,1)\\
0 & \mathrm{\ otherwise}%
\end{array}
\right.
\]
Invoking (\ref{discupdate}) for $x(t,\theta)$ gives
\begin{equation}
x^{\prime}(\tau_{k}^{+},\theta)=x^{\prime}(\tau_{k}^{-},\theta)+\Delta
f_{x}(q_{k-1},q_{k},\theta){\tau}_{k}^{\prime}.\label{xdiscupdate}\end{equation}
By (\ref{finaltkprime}), we only need to take care of the
transitions caused
by $E_{2}$ and $E_{3}$ as in other cases $\Delta f_{x}(q_{k-1},q_{k}%
,\theta)\tau_{k}^{\prime}=0$. Since neither $E_{2}$ nor $E_{3}$ appear in the
transitions with resets, they cannot create a chain of simultaneous
transitions, thereby leaving us with transition $(1,2)$ for $E_{2}$ and
$(1,0)$ for $E_{3}$. In the first case, we get $\Delta f_{x}(1,2,\theta
)\tau_{k}^{\prime}=1-x^{\prime}(\tau_{k}^{-},\theta)$ and in the latter case
we get $\Delta f_{x}(1,0,\theta)\tau_{k}^{\prime}=-x^{\prime}(\tau_{k}%
^{-},\theta)$. Inserting these results in (\ref{xdiscupdate}) yields
\begin{equation}
x^{\prime}(\tau_{k}^{+},\theta)=\left\{
\begin{array}
[c]{ll}%
1 & \mathrm{\ if\ }(q_{k-1},q_{k})=(1,2)\\
0 & \mathrm{\ if\ }(q_{k-1},q_{k})=(1,0)\\
x^{\prime}(\tau_{k}^{-},\theta) & \mathrm{\ otherwise}%
\end{array}
\right.\label{examplediscreteupdate}
\end{equation}

There is no need to consider (\ref{contupdate}) in this case,
since $\frac{d\mathrm{\mathbf{f}}_{q}}{d\theta}=0$ for all $q$.
Therefore, we are in the position to fully evaluate the sample
derivative estimates (\ref{Qprime}) and \eqref{Lprime}.

By (\ref{examplediscreteupdate}), after the buffer becomes empty (transition
(1,0) through event $E_{3}$), $x^{\prime}(t,\theta)$ becomes and stays at 0
until a transition $(1,2)$ occurs through $E_{2}$. If this happens,
$x^{\prime}(t,\theta)$ resets to 1 in (\ref{examplediscreteupdate}) and
remains constant until the buffer becomes empty again. Therefore, we need only
consider those nonempty periods $[\xi_{n},\eta_{n})$ in which a transition
$(1,2)$ occurs. If this happens, we calculate the length of the interval
between the first such transition until the next time the buffer becomes
empty. $Q^{\prime}(T,\theta)$ in (\ref{Qprime}) is the sum of lengths of these
intervals, i.e.,
\begin{equation}
Q^{\prime}(T,\theta)=\frac{1}{T}\sum_{n=1}^{N}\mathbf{1}_{FP}(n)(\eta_{n}%
-\nu_{n,1}) \label{exampleQprime}%
\end{equation}
where $\mathbf{1}_{FP}(n)=1$ if there exists a transition $(1,2)$ in the
non-empty period $[\xi_{n},\eta_{n})$ and 0 otherwise.

Next, to evaluate \eqref{Lprime}, notice that at $t=\sigma_{n,m}$
(end of stay at $q=2$ in Fig. \ref{fig:singleclassSFM}) a
transition to $q=1$ can occur in two ways: (a) Through $E_{4}$ or
$E_{5}$ (transition (2,2)) and violating the invariant condition
$[\alpha>\beta]$ which immediately fires transition $(2,1)$; (b)
Directly, by $E_{1}$ (transition (2,1)). These three possibilities
are associated with zeros in (\ref{finaltkprime}), so we have
$\sigma_{n,m}^{\prime}=0$. Regarding the term
$-\nu_{n,m}^{\prime}\ell(\nu_{n,m}^{+},\theta)$ in \eqref{Lprime},
the event at $\nu_{n,m}$ is $E_{2}$. By (\ref{finaltkprime}), we
have
$\nu_{n,m}^{\prime}=-\frac{x^{\prime}(\nu_{n,m}^{-},\theta)-1}{\alpha
(\nu_{n,m}^{-})-\beta(\nu_{n,m}^{-})}$. Since by (\ref{lossrate}),
$\ell
(\nu_{n,m}^{+},\theta)=\alpha(\nu_{n,m}^{+})-\beta(\nu_{n,m}^{+})$
and by
Assumption \ref{assump:nosimul}, $\alpha(\nu_{n,m}^{+})-\beta(\nu_{n,m}%
^{+})=\alpha(\nu_{n,m}^{-})-\beta(\nu_{n,m}^{-})$, we find that $-\nu
_{n,m}^{\prime}\ell(\nu_{n,m}^{+},\theta)=x^{\prime}(\nu_{n,m}^{-},\theta)-1$.
We have already shown in (\ref{examplediscreteupdate}) that in a non-empty
period $[\xi_{n},\eta_{n})$, $x^{\prime}(t,\theta)=0$ for all $t\in\lbrack
\xi_{n},\nu_{n,1})$ and $x^{\prime}(t,\theta)=1$, $t\in\lbrack\nu_{n,1}%
,\eta_{n})$. Hence, $x^{\prime}(\nu_{n,m}^{-},\theta)=0$ when
$m=1$ and $x^{\prime}(\nu_{n,m}^{-},\theta)=1$, otherwise.
Combining all results into \eqref{Lprime}, we find that
\[
L^{\prime}(T,\theta)=\sum_{n=1}^{N}\sum_{m=1}^{M_{n}}-\nu_{n,m}^{\prime}%
\ell(\nu_{n,m}^{+},\theta)=-N_{F}%
\]
where $N_{F}$ is the number of non-empty intervals with at least
one full period. These results recover those in \citep[pp.
700--703]{cassbook} and \citep{CassWarMelSunPan02}. Note that
$Q^{\prime}(T,\theta)$ and $L^{\prime}(T,\theta)$ are independent
of $f_{\alpha}$ and $f_{\beta}$, i.e., these sensitivity estimates
are independent of the random arrival and service processes, a
fundamental robustness property of IPA.

\section{Conclusions}

We have introduced a general framework suitable for analysis and
on-line optimization of Stochastic Hybrid Systems (SHS) which
facilitates the use of Infinitesimal Perturbation Analysis (IPA).
In doing so, we modified the previous structure of a Stochastic
Hybrid Automaton (SHA) and showed that every transition is
associated with an explicit event which is defined through a guard
function. This also enables us to uniformly treat all events
observed on the sample path of the SHS and makes it possible to
develop a unifying matrix notation for IPA equations which
eliminates the need for the case-by-case analysis based on event
classes as in prior work involving IPA for SHS.
\bibliographystyle{ifacconf}
\bibliography{AliBIB}

\begin{thebibliography}{13}
\providecommand{\natexlab}[1]{#1}
\providecommand{\url}[1]{\texttt{#1}}
\expandafter\ifx\csname urlstyle\endcsname\relax
  \providecommand{\doi}[1]{doi: #1}\else
  \providecommand{\doi}{doi: \begingroup \urlstyle{rm}\Url}\fi

\bibitem[{Bujorianu} and {Lygeros}(2006)]{BujLyg06}
L.~M. {Bujorianu} and J.~{Lygeros}.
\newblock Toward a general theory of stochastic hybrid systems.
\newblock In \emph{Stochastic Hybrid Systems: Theory and Safety Critical
  Applications}, volume 337 of \emph{Lect. Notes in Cont. and Inf. Sci.}, pages
  3--30. Springer Verlag, Berlin, July 2006.

\bibitem[Bujorianu and Lygeros(2004)]{BujLyg04}
M.L. Bujorianu and J.~Lygeros.
\newblock General stochastic hybrid systems: modelling and optimal control.
\newblock In \emph{43rd IEEE Conf. on Dec. and Cont.}, volume~2, pages 1872 --
  1877 Vol.2, dec. 2004.

\bibitem[Cassandras and Lafortune(2006)]{cassbook}
C.~G. Cassandras and Stephane Lafortune.
\newblock \emph{Introduction to Discrete Event Systems}.
\newblock Springer-Verlag New York, Inc., Secaucus, NJ, USA, 2006.

\bibitem[Cassandras and Lygeros(2006)]{SHSbook}
C.~G. Cassandras and J.~Lygeros.
\newblock \emph{Stochastic Hybrid Systems}.
\newblock Taylor and Francis, 2006.

\bibitem[Cassandras et~al.(2002)Cassandras, Wardi, Melamed, Sun, and
  Panayiotou]{CassWarMelSunPan02}
C.~G. Cassandras, Y.~Wardi, B.~Melamed, G.~Sun, and C.~G. Panayiotou.
\newblock Perturbation analysis for on-line control and optimization of
  stochastic fluid models.
\newblock \emph{IEEE Trans. on Aut. Cont.}, AC-47\penalty0 (8):\penalty0
  1234--1248, 2002.

\bibitem[Cassandras et~al.(2009)Cassandras, Wardi, Panayiotou, and
  Yao]{CassWarPanYao09}
C.~G. Cassandras, Y.~Wardi, C.~G. Panayiotou, and C.~Yao.
\newblock Perturbation analysis and optimization of stochastic hybrid systems.
\newblock \emph{European J. of Cont.}, 16\penalty0 (6):\penalty0 642--664,
  2009.

\bibitem[Ghosh et~al.(1993)Ghosh, Arapostathis, and Marcus]{GhoshArapMarc93}
M.~K. Ghosh, A.~Arapostathis, and S.~I. Marcus.
\newblock Optimal control of switching diffusions with application to flexible
  manufacturing systems.
\newblock \emph{SIAM}, 31\penalty0 (5):\penalty0 1183--1204, 1993.

\bibitem[Ghosh et~al.(1997)Ghosh, Arapostathis, and Marcus]{GhoshArapMarc97}
M.~K. Ghosh, A.~Arapostathis, and S.~I. Marcus.
\newblock Ergodic control of switching diffusions.
\newblock \emph{SIAM}, 35\penalty0 (6):\penalty0 1952--1988, 1997.

\bibitem[Glynn(1989)]{glynn89}
P.W. Glynn.
\newblock A gsmp formalism for discrete event systems.
\newblock \emph{Proc. of the IEEE}, 77\penalty0 (1):\penalty0 14 --23, jan
  1989.

\bibitem[Hespanha(2004)]{Hespanha04}
J.~P. Hespanha.
\newblock Stochastic hybrid systems: Application to communication networks.
\newblock In \emph{Hybrid Systems: Comput. and Cont.}, pages 387--401.
  Springer-Verlag, 2004.

\bibitem[Koutsoukos(2005)]{Koutsoukos05}
X.D. Koutsoukos.
\newblock Optimal control of stochastic hybrid systems based on locally
  consistent markov decision processes.
\newblock In \emph{Proc. of the IEEE Int'l Symp. on Intel. Cont.}, pages 435
  --440, 2005.

\bibitem[Simi\'{c} et~al.(2000)Simi\'{c}, Johansson, Sastry, and
  Lygeros]{Sim00}
S.~Simi\'{c}, K.~Johansson, S.~Sastry, and J.~Lygeros.
\newblock Towards a geometric theory of hybrid systems.
\newblock In \emph{Hybrid Systems: Comp. and Cont.}, volume 1790 of \emph{Lect.
  Notes in Computer Sci.}, pages 421--436. Springer, 2000.

\bibitem[Yao and Cassandras(2011)]{YaoCass11}
C.~Yao and C.~G. Cassandras.
\newblock Perturbation analysis of stochastic hybrid systems and applications
  to resource contention games.
\newblock \emph{Frontiers of Electrical and Electronic Eng. in China},
  6\penalty0 (3):\penalty0 453--467, 2011.

\end{thebibliography}

\end{document}